\documentclass[10pt]{article}

\usepackage[english]{babel}
\usepackage{amsfonts,amsmath,amssymb}
\usepackage{geometry}

\newtheorem{proposition}{Proposition}
\newtheorem{corollary}{Corollary}
\newtheorem{remark}{Remark}
\newtheorem{definition}{Definition}

\newtheorem{lemma}{Lemma}
\newtheorem{example}{Example}
\newenvironment{proof}[1][Proof:]{\begin{trivlist}
\item[\hskip \labelsep {\bfseries #1}]}{\end{trivlist}}

\newcommand{\id}{{\mathbf 1}}
\newcommand{\Hom}{\mathrm{Hom}}
\newcommand{\Aut}{\mathrm{Aut}}
\newcommand{\ad}{\operatorname{ad}}
\newcommand{\End}{\mathrm{End}}
\newcommand{\Der}{\mathrm{Der}}
\newcommand{\F}{\mathbb{F}}
\newcommand{\codim}{\operatorname{codim}}
\newcommand{\onto}{{\!\mapsto\!\!\!\!\!\rightarrow\!}}
\newcommand{\LO}{\mathrm{LO}}

\begin{document}

\bibliographystyle{plain}

\setcounter{page}{1}

\thispagestyle{empty}

\title{The structure of almost Abelian Lie algebras}

\author{
Zhirayr Avetisyan
}

\maketitle

\begin{abstract}
An almost Abelian Lie algebra is a non-Abelian Lie algebra with a codimension 1 Abelian ideal. Most 3-dimensional real Lie algebras are almost Abelian, and they appear in every branch of physics that deals with anisotropic media - cosmology, crystallography etc. In differential geometry and theoretical physics, almost Abelian Lie groups have given rise to some of the simplest solvmanifolds on which various geometric structures such as symplectic, K\"ahler, spin etc., are currently studied in explicit terms. However, a systematic study of almost Abelian Lie groups and algebras from mathematics perspective has not been carried out yet, and the present paper is the first step in addressing this wide and diverse class of groups and algebras.

The present paper studies the structure and important algebraic properties of almost Abelian Lie algebras of arbitrary dimension over any field of scalars. A classification of almost Abelian Lie algebras is given. All Lie subalgebras and ideals, automorphisms and derivations, Lie orthogonal operators and quadratic Casimir elements are described exactly.
\end{abstract}


\section{Introduction}

In its most common definition an almost Abelian Lie algebra is a non-Abelian Lie algebra $\mathbf{L}$ over a field $\F$ which has a codimension 1 Abelian ideal. However, it was shown in \cite{BuCe12} (Proposition 3.1) that if a Lie algebra has a codimension 1 Abelian subalgebra then it necessarily also has a codimension 1 Abelian ideal (which may be different from the original subalgebra). The result is announced for fields $\F$ of characteristic zero and in finite dimensions, but in fact the proof works with minor adaptations in the most general case. Thus, the seemingly weaker definition of an almost Abelian Lie algebra as possessing a codimension 1 Abelian subalgebra is in fact equivalent to the original one. For purely aesthetic purposes in this paper we will adopt the latter, apparently weaker but factually equivalent definition.

There is certain controversy regarding the notion of an almost Abelian Lie algebra in the literature. Some authors (e.g., \cite{HLP90}, \cite{GlKo69}, \cite{Kol65}) define an almost Abelian Lie algebra as a Lie algebra for which there exists a basis $\{e_0,e_i\}$, $i=1,...,n$ such that $[e_0,e_i]=e_i$ and $[e_i,e_j]=0$ for all $i$, i.e., $\ad_{e_0}=\id$. This definition is too strict in that it misses all other possibilities for the operator $\ad_{e_0}$ not similar to identity. On the other hand, some sources prefer to include Abelian Lie algebras in the class of almost Abelian Lie algebras (e.g., \cite{Pop13}), but we prefer to not do so because of principal structural differences between the two.

Speaking of applications, let us consider only real and complex Lie algebras in this paragraph. In low dimensions almost Abelian Lie algebras are vastly represented. The only 2-dimensional non-Abelian Lie algebra is almost Abelian, while 6 out of 9 classes (after Bianchi) of 3-dimensional real Lie algebras are almost Abelian. On the other hand, since most physical systems are $n=1,2,3$ dimensional, in absence of rotational symmetries a homogeneous (anisotropic) system is described by a low dimensional Lie algebra. It is thus not a surprise that almost Abelian Lie algebras are widely used in cosmology, where they represent the symmetries of the universe at large scale (\cite{ElEl98}, \cite{Osi73}, \cite{Pet59}, \cite{Rya75}, \cite{AvVe13} and many others), or crystallography, where they model the symmetries of an ideal solid (\cite{Par16} and references therein). As far as applications in pure mathematics are concerned, one particular almost Abelian Lie algebra together with its Lie groups is distinguished - the 3-dimensional Heisenberg algebra (higher dimensional Heisenberg algebras are not almost Abelian). Thorough studies of the Heisenberg algebra and the corresponding Lie group can be found, for instance, in \cite{Fol89} and \cite{Tha98}. Taking roots in the foundations of quantum mechanics, this algebra has become the classical setting for non-commutative analysis. We refer to \cite{FiRu16} for a thorough study of quantization and pseudodifferential calculus on the Heisenberg group (among other nilpotent groups). It is only natural to try and extend these beautiful results to general almost Abelian groups, but that should wait until a comprehensive study of almost Abelian Lie algebras and groups is available.

Higher dimensional almost Abelian Lie groups and algebras have gained in popularity in the last two decades, with at least a dozen papers dealing with the subject written in the last two years only. One context of interest is compact solvmanifolds. A solvmanifold is a homogeneous space $G/N$ with $G$ a simply connected solvable Lie group and $N\subset G$ a discrete subgroup. Almost Abelian groups $G$ are special in that, together with nilpotent groups, these are the only solvable Lie groups for which there is a practically useful necessary and sufficient condition for the solvmanifold $G/N$ to be compact \cite{Boc16}. More generally, almost Abelian groups and algebras are unique in their explicit tractability combined with diversity of properties they can possess. A plethora of work in differential geometry and theoretical physics has been devoted to various geometrical constructions on almost Abelian solvmanifolds such as symplectic, K\"ahler, spin, $\mathrm{G}_2$ or $\mathrm{SU}(3)$ structures, various flows etc. \cite{Fre12}, \cite{AnOr17}, \cite{AGMP11}, \cite{CoMa12}, \cite{LaWi19},\cite{FSW19}, \cite{Par21}, \cite{Sta20}, \cite{FiPa20a}, \cite{FiPa20}, \cite{BDV19}, \cite{BaFi18}, \cite{FrSw18}, \cite{BDV18}. In spite of this wide spectrum of interest and applications, there exists to the date no comprehensive study of almost Abelian Lie algebras and groups in the literature. The present paper is the first step in the systematic study of almost Abelian Lie algebras, groups and their homogeneous spaces.

The simple idea behind the structure of an almost Abelian Lie algebra has equivalent manifestations in several mathematical and physical formalisms. For instance, it can be given the following interpretation in the theory of integrable systems and PDE. Let $(C^\infty(M),\{,\})$ be a Poisson algebra over an $n$-dimensional manifold $M$, and let $H\in C^\infty(M)$ be a Hamiltonian. The Hamiltonian system evolving according to $H$ is Liouville integrable if $H$ is contained in an Abelian subalgebra of $(C^\infty(M),\{,\})$ (considered as a Lie algebra) generated by at least $n$ functionally independent integrals of motion. In particular, if $H$ is invariant under the action of a Lie algebra $\mathbf{L}$ of Noether symmetries which has an at least $n-1$-dimensional Abelian subalgebra then the system is integrable in Noether integrals. Thus, a (left) invariant Hamiltonian system on a Lie group is integrable in Noether symmetries if and only if the group is either Abelian or almost Abelian (this is one occasion where including Abelian Lie algebras in almost Abelian ones would make the statement more elegant, but we will resist the temptation).

Another way an almost Abelian Lie algebra can be viewed is a linear dynamical system. Let $\F$ be a field, and $\mathbf{V}$ an $\F$-vector space. A discrete linear dynamical system is an equation of the form
$$
v(n+1)=\mbox{T}v(n),\quad v:\mathbb{N}_0\mapsto\mathbf{V},
$$
where $\mbox{T}\in\End_\F(\mathbf{V})$ is a linear operator. To give our dynamical system a Hamiltonian flavour consider the extended $\F$-vector space
$$
\mathbf{L}=\mathbf{V}\oplus_\F\F H
$$
with a distinguished element $H\in\mathbf{L}$ called the {\it Hamiltonian}. Make $\mathbf{L}$ a Lie algebra by setting
$$
[u,v]=0,\quad [H,v]=\ad_Hv=\mbox{T}v,\quad\forall v,u\in\mathbf{V}.
$$
Now if $\mbox{T}\neq0$ then $\mathbf{L}$ is an almost Abelian Lie algebra that encodes the dynamics of the original system (if $\mbox{T}=0$ then $\mathbf{L}$ is Abelian and the dynamics is trivial),
$$
v(n+1)=[H,v(n)].
$$
If the vector space $\mathbf{V}$ is endowed with a topology then the continuous counterpart of the system can be considered as well,
$$
\frac{\operatorname{d} v(t)}{\operatorname{d} t}=[H,v(t)].
$$
Thus, an almost Abelian Lie algebra is also a stationary linear dynamical system, and the study of its algebraic properties yields to the understanding of the dynamics. Linear dynamical systems are well understood and solved in terms of matrix algebra. But things are not as trivial when $\mathbf{V}$ is infinite dimensional, or when $\F$ is neither $\mathbb{R}$ nor $\mathbb{C}$. Matrix tools such as determinants, traces and transpositions need not exist. Moreover, even advanced methods like canonical forms (e.g., Jordan or Frobenius) become non-trivial issues and need not exist. No reference can be made to the spectral theory of Hilbert space operators because $\mathbf{V}$ is merely a vector space. And even in the simplest case of a real or complex finite dimensional diagonalizable matrix $\mbox{T}$ many algebraic aspects need to be described explicitly and in a systematic manner.

The smallest (in terms of dimensions) almost Abelian Lie algebra is the only (up to isomorphism) non-Abelian Lie algebra of dimension 2, which is the Lie algebra $\mathbf{ax\!+\!b}_\F$ of the group of affine transformations on the $\F$-line. This algebra can be described in a basis $\{e_0,e_1\}$ with relations $[e_0,e_1]=e_1$, i.e., $\ad_{e_0}=\id$. In dimension 3 the almost Abelian Lie algebras are the non-Abelian solvable Bianchi algebras $\mbox{Bi}(\mbox{II})$-$\mbox{Bi}(\mbox{VII})$ ($\mbox{Bi}(\mbox{I})$ is Abelian and $\mbox{Bi}(\mbox{VIII})$ and $\mbox{Bi}(\mbox{IX})$ are simple, so we will not consider them here) which were first classified for $\F=\mathbb{R}$ by Bianchi \cite{Bia98} (see \cite{GKKL14} for a more modern and elegant approach). The classification corresponds to the similarity classes of $2\times2$ matrices $\ad_{e_0}$. The most prominent of them is the nilpotent Heisenberg algebra $\mathbf{H}_\F$ which arises from the canonical commutation relations in quantum mechanics,
$$
e_0=\frac\partial{\partial x},\quad e_1=\id,\quad e_2=x,\quad[e_0,e_2]=e_1,\quad[e_0,e_1]=[e_1,e_2]=0.
$$

In this paper we consider almost Abelian Lie algebras $\mathbf{L}$ of any dimension $\dim_\F\mathbf{L}=1+\aleph$ over any field $\F$, so that the codimension 1 Abelian ideal has dimension $\aleph$. Note that we identify index sets with their cardinalities, therefore we use the symbol $\aleph$ to denote both the cardinality of a Hamel basis in a vector space and an actual index set of that cardinality. This ultimate generality enforces strict limitations on the arsenal of algebraic methodologies that can be utilized. No integers other than $0$ and $\pm1$ are assumed to be contained in $\F$ a priory in order to avoid a potential conflict with the characteristic of $\F$. No reference can be made to polynomial algebraic and spectral methods such as eigenvalues or factorizations so that algebraic and topological aspects of $\F$ such as algebraic extensions or topological closedness are irrelevant. And of course, no dimensional or matrix arguments such as rank-nullity or canonical forms can be used in this generality. This may have made some otherwise simple questions unexpectedly tricky. Some structures such as homomorphisms become somewhat richer in infinite dimensions.

The paper is structured as follows. First in section \ref{LieAlg} we introduce notations by reminding some notions and facts from the elementary Lie theory which are used in the body. This makes the paper essentially self-contained with little reference to the standard background in the field, and thereby more accessible to non-specialists (especially non-mathematicians) who may turn out to be a majority among those interested in the subject. Then in section \ref{AALieAlg} we define almost Abelian Lie algebras, discuss their subalgebras and ideals, and decompose them into an indecomposable core and a central extension. Further in section \ref{HomClass} we describe homomorphisms between and automorphisms of almost Abelian Lie algebras which then naturally leads to a classification up to isomorphism. Afterwards in section \ref{DerOrth} we describe derivations and Lie orthogonal operators on almost Abelian Lie algebras. Finally in section \ref{Casimir} we study the centre of the universal enveloping algebra of an almost Abelian Lie algebra.


\section{Lie algebras\label{LieAlg}}

Lie algebras will be the main subject of this exposition. For a thorough exposition of Lie algebras and Lie groups we refer to the all-time classical monograph by Knapp \cite{Kna02}. Here we will collect some definitions and notations to be used in the sequel that may not be completely obvious.

We will consider Lie algebras of any dimension over an arbitrary field $\F$. We will write $\mathbf{L}_1\simeq\mathbf{L}_2$ if two Lie algebras $\mathbf{L}_1$ and $\mathbf{L}_2$ are isomorphic. The vector space direct sum $\mathbf{L}_1\oplus_\F\mathbf{L}_2$ will be distinguished from the Lie algebra direct sum $\mathbf{L}_1\oplus\mathbf{L}_2$. The semidirect product (or sum) of two Lie algebras $\mathbf{L}_1$ and $\mathbf{L}_2$ will be denoted by $\mathbf{L}_1\rtimes\mathbf{L}_2$.

A Lie algebra $\mathbf{L}$ is called decomposable if it decomposes into a direct sum $\mathbf{L}=\mathbf{L}_1\oplus\mathbf{L}_2$ of two commuting non-trivial Lie subalgebras $\mathbf{L}_1\neq0$ and $\mathbf{L}_2\neq0$. Otherwise $\mathbf{L}$ is called indecomposable.

The nilradical $\mathfrak{nil}(\mathbf{L})\subset\mathbf{L}$ of a Lie algebra $\mathbf{L}$ is a maximal nilpotent ideal. If it exists then it is unique, because the sum of two nilpotent ideals is again a nilpotent ideal. It always exists when $\dim_\F\mathbf{L}<\infty$.

Automorphisms of a given Lie algebra comprise a group denoted by $\Aut(\mathbf{L})$.

The set $\Der(\mathbf{L})\subset\End_\F(\mathbf{L})$ of all derivations of a given Lie algebra $\mathbf{L}$ comprises a Lie subalgebra of $\End_\F(\mathbf{L})$ under commutation. For every $X\in\mathbf{L}$, the operator $\ad_X\in\End_\F(\mathbf{L})$ is a derivation. Thus, the image of the adjoint representation $\ad_\mathbf{L}\subset\Der(\mathbf{L})$ is a Lie subalgebra called the algebra of inner derivations. Those derivations $D\in\Der(\mathbf{L})\setminus\ad_\mathbf{L}$ are called outer derivations.

A linear operator $\mbox{T}\in\End_\F(\mathbf{L})$ on a Lie algebra $\mathbf{L}$ is called Lie orthogonal (see \cite{Pop13}) if
$$
[\mbox{T}X,\mbox{T}Y]=[X,Y],\quad\forall X,Y\in\mathbf{L}.
$$
Obviously the set $\LO(\mathbf{L})\subset\End_\F(\mathbf{L})$ of all Lie orthogonal operators on a given Lie algebra $\mathbf{L}$ is a monoid under composition with the identity map being the unit.

For a Lie algebra $\mathbf{L}$ consider the tensor algebra $\operatorname{T}(\mathbf{L})$ which is a unital associative graded non-commutative $\F$-algebra, and consider the ideal $\mathbf{I}_\mathbf{L}\subset\operatorname{T}(\mathbf{L})$ generated by elements of the form
$$
X\otimes Y-Y\otimes X-[X,Y]\in\operatorname{T}(\mathbf{L}),\quad\forall X,Y\in\mathbf{L}.
$$
The quotient algebra
$$
\operatorname{U}(\mathbf{L})\doteq\operatorname{T}(\mathbf{L})/\mathbf{I}_\mathbf{L}
$$
is called the universal enveloping algebra of $\mathbf{L}$. The natural inclusion $\mathbf{L}\hookrightarrow\operatorname{T}(\mathbf{L})$ gives rise to an injective map $\mathfrak{h}:\mathbf{L}\to\operatorname{U}(\mathbf{L})$ such that
$$
\mathfrak{h}(X)\mathfrak{h}(Y)-\mathfrak{h}(Y)\mathfrak{h}(X)-\mathfrak{h}([X,Y])=0,\quad\forall X,Y\in\mathbf{L}.
$$
Denote by
$$
\mathcal{Z}(\operatorname{U}(\mathbf{L}))=\left\{x\in\operatorname{U}(\mathbf{L})|\quad xy=yx,\quad\forall y\in\operatorname{U}(\mathbf{L})\right\}
$$
the centre of the universal enveloping algebra.

Let $D\in\Der(\mathbf{L})$ be a derivation, and let $\operatorname{U}(\mathbf{L})$ be the universal enveloping algebra. Using the map $\mathfrak{h}:\mathbf{L}\to\operatorname{U}(\mathbf{L})$ we can consider $\mathbf{L}$ as a subspace of $\operatorname{U}(\mathbf{L})$. We can extend the action of $D$ to $\operatorname{U}(\mathbf{L})$ by requiring it to be a derivation,
$$
D(x*y)=Dx*y+x*Dy,\quad\forall x,y\in\operatorname{U}(\mathbf{L}).
$$
In doing so we start from extending $D$ to entire $\operatorname{T}(\mathbf{L})$ as a derivation and using the fact that $D\mathbf{I}_\mathbf{L}\subset\mathbf{I}_\mathbf{L}$ which follows from
$$
D(X\otimes Y-Y\otimes X-[X,Y])=DX\otimes Y-Y\otimes DX-[DX,Y]+X\otimes DY-DY\otimes X-[X,DY].
$$
Thus, $\Der(\mathbf{L})\hookrightarrow\Der(\operatorname{U}(\mathbf{L}))$.


\section{The structure of almost Abelian Lie algebras\label{AALieAlg}}

In this section we will start serving our main purpose by studying the structure of almost Abelian Lie algebras. Abelian Lie algebras are simple objects characterized by their underlying vector space only, and are therefore often denoted by $\mathbf{L}\doteq\F^\aleph$ with $\aleph=\dim_\F\mathbf{L}$ (the cardinality of an $\F$-basis). Arguably the next simplest Lie algebras are almost Abelian Lie algebras.

\begin{definition} A non-Abelian Lie algebra $\mathbf{L}$ is called almost Abelian if it contains a codimension 1 Abelian subalgebra.
\end{definition}
This implies that $\mathbf{L}=\F^\aleph\oplus_\F\F$ as $\F$-vector spaces where $\aleph=\dim_\F\mathbf{L}-1$ and $[\mathbf{L},\F^\aleph]\neq0$. The following is a minimal modification of Proposition 3.1 in \cite{BuCe12}.

\begin{proposition} An almost Abelian Lie algebra $\mathbf{L}$ over $\F$ has a codimension 1 Abelian ideal, and is therefore isomorphic to the semidirect product
$$
\mathbf{L}\simeq\F^\aleph\rtimes\F,\quad\aleph=\dim_\F\mathbf{L}-1.
$$
\end{proposition}
\begin{proof} The proof of the existence of a codimension 1 Abelian ideal $\F^\aleph$ is a word-by-word adaptation of that in the above mentioned reference to the infinite dimensional setting. $\mathbf{L}$ is a semidirect product because the quotient Lie algebra $\mathbf{L}/\F^\aleph\simeq\F$ is isomorphic to the subalgebra $0\oplus\F\subset\mathbf{L}$.
\end{proof}

Henceforth we will assume $\mathbf{L}=\mathbf{V}\rtimes\F e_0$ for an element $e_0\in\mathbf{L}$ and the vector space $\mathbf{V}=\F^\aleph$, and will choose an $\F$-basis $\{e_i\}_{i\in\aleph}$ of $\mathbf{V}$. By Proposition \ref{2cod1AbIdProp} that will come later the codimension 1 Abelian ideal $\mathbf{V}$ is essentially unique and thus, (with one notable exception that we will always bear in mind) no loss of generality occurs in writing $\mathbf{L}=\mathbf{V}\rtimes\F e_0$.  We will write every $X\in\mathbf{L}$ as $X=(v,t)\in\mathbf{V}\oplus_\F\F$. The Lie algebra structure of $\mathbf{L}$ is completely determined by
$$
[e_0,v]=\ad_{e_0}v,\quad v\in\mathbf{V}.
$$
Here $\ad_{e_0}\in\End_\F(\mathbf{V})$ may be any nonzero linear operator on $\mathbf{V}$.

\begin{remark}\label{VadPairRemark} Thus, an almost Abelian Lie algebra is given by a pair $(\mathbf{V},\ad_{e_0})$ of a vector space $\mathbf{V}$ and a nonzero linear operator $\ad_{e_0}\in\End_\F(\mathbf{V})\setminus\{0\}$. Every such pair gives rise to an almost Abelian Lie algebra, but different pairs can yield isomorphic Lie algebras.
\end{remark}
Denote by $\mathcal{Z}(\mathbf{L})$ the centre of the Lie algebra $\mathbf{L}$,
$$
\mathcal{Z}(\mathbf{L})\doteq\left\{X\in\mathbf{L}|\quad[X,\mathbf{L}]=0\right\}.
$$

\begin{remark}\label{ZRemark} For an almost Abelian Lie algebra $\mathbf{L}=\mathbf{V}\rtimes\F e_0$ we have
$$
\mathcal{Z}(\mathbf{L})=\ker\ad_{e_0},\quad[\mathbf{L},\mathbf{L}]=\ad_{e_0}\mathbf{V},\quad[[\mathbf{L},\mathbf{L}],[\mathbf{L},\mathbf{L}]]=0,
$$
$$
[\mathbf{L},...,[\mathbf{L},\mathbf{L}]...]=\ad_{e_0}^n\mathbf{V},\quad n\in\mathbb{N}.
$$
\end{remark}
Thus, an almost Abelian $\mathbf{L}$ is always 2-step solvable, and it is nilpotent if and only if the operator $\ad_{e_0}$ is. The nilradical is $\mathfrak{nil}(\mathbf{L})=\mathbf{L}$ if $\mathbf{L}$ is nilpotent and $\mathfrak{nil}(\mathbf{L})=\mathbf{V}$ otherwise. In the latter case no ambiguity arises because as we will see later (Proposition \ref{2cod1AbIdProp}), in the rare case when there are more than one codimension 1 Abelian ideals $\mathbf{V}$, $\mathbf{L}$ is necessarily nilpotent. The adjoint representation is given by
\begin{equation}
\ad_X=\begin{pmatrix}
t\ad_{e_0} & -\ad_{e_0}v\\
0 & 0
\end{pmatrix},\quad\forall X=(v,t)\in\mathbf{L}.\label{adRep}
\end{equation}
We observe immediately that if $\ker\ad_{e_0}=\mathcal{Z}(\mathbf{L})\neq0$ then $\ad$ is not faithful.

\begin{proposition} Every almost Abelian Lie algebra $\mathbf{L}=\mathbf{V}\rtimes\F e_0$ has a faithful representation in $\End_\F(\F\oplus_\F\mathbf{V})$ given by
$$
X=(v,t)\mapsto\begin{pmatrix}
0 & 0\\
v & t\ad_{e_0}
\end{pmatrix},\quad\forall X\in\mathbf{L}.
$$
\end{proposition}
\begin{proof} We first observe that this map is $\F$-linear and injective, and that
$$
\begin{pmatrix}
0 & 0\\
\mathbf{V} & 0
\end{pmatrix}\simeq\mathbf{V},\quad\begin{pmatrix}
0 & 0\\
0 & \F\ad_{e_0}
\end{pmatrix}\simeq\F e_0
$$
as Abelian Lie algebras. Then we check that
$$
\ad_{e_0}v=[e_0,v]=\left[\begin{pmatrix}
0 & 0\\
0 & \ad_{e_0}
\end{pmatrix},\begin{pmatrix}
0 & 0\\
v & 0
\end{pmatrix}\right]=\begin{pmatrix}
0 & 0\\
\ad_{e_0}v & 0
\end{pmatrix}
$$
which proves the statement.
\end{proof}
The next proposition describes all possible Lie subalgebras and ideals of a given almost Abelian Lie algebra.

\begin{proposition}\label{SubalgIdealProp} Every Lie subalgebra $\mathbf{L}_0\subset\mathbf{L}$ of an almost Abelian Lie algebra $\mathbf{L}=\mathbf{V}\rtimes\F e_0$ has one of the following forms:
\begin{itemize}
\item An Abelian Lie subalgebra $\mathbf{L}_0=\mathbf{W}$ with any subspace $\mathbf{W}\subset\mathbf{V}$

\item An Abelian Lie subalgebra $\mathbf{L}_0=\mathbf{W}\oplus\F e_1$ with any $e_1\in e_0+\mathbf{V}$ and any subspace $\mathbf{W}\subset\ker\ad_{e_0}$

\item An almost Abelian Lie subalgebra $\mathbf{L}_0=\mathbf{W}\rtimes\F e_1$ with any $e_1\in e_0+\mathbf{V}$ and any $\ad_{e_0}$-invariant subspace $\mathbf{W}\subset\mathbf{V}$ with $\mathbf{W}\not\subset\ker\ad_{e_0}$
\end{itemize}
Every ideal $\mathbf{I}\subset\mathbf{L}$ has one of the following forms:
\begin{itemize}
\item An Abelian ideal $\mathbf{I}=\mathbf{W}$ with any $\ad_{e_0}$-invariant subspace $\mathbf{W}\subset\mathbf{V}$

\item An Abelian ideal $\mathbf{I}=\mathbf{W}\oplus\F e_1$ with any $e_1\in e_0+\mathbf{V}$ and any subspace $\mathbf{W}\subset\ker\ad_{e_0}$ with $\ad_{e_0}\mathbf{V}\subset\mathbf{W}$

\item An almost Abelian ideal $\mathbf{I}=\mathbf{W}\rtimes\F e_1$ with any $e_1\in e_0+\mathbf{V}$ and any subspace $\mathbf{W}\subset\mathbf{V}$ with $\ad_{e_0}\mathbf{V}\subset\mathbf{W}\not\subset\ker\ad_{e_0}$
\end{itemize}
\end{proposition}
\begin{proof} Let $\mathbf{L}_0\subset\mathbf{L}$ be an arbitrary subspace. Then either $\mathbf{L}_0=\mathbf{W}\subset\mathbf{V}$ or $\mathbf{L}_0=\mathbf{W}\oplus_\F\F e_1$ for some $e_1\in e_0+\mathbf{V}$ and a possibly trivial subspace $\mathbf{W}\subset\mathbf{V}$. In the former case $\mathbf{L}_0$ is readily an Abelian Lie subalgebra. In the latter case, in order for $\mathbf{L}_0$ to be a Lie subalgebra we need
$$
[\mathbf{W}\oplus_\F\F e_1,\mathbf{W}\oplus_\F\F e_1]=\ad_{e_1}\mathbf{W}=\ad_{e_0}\mathbf{W}\subset\mathbf{W}.
$$
If $\mathbf{W}\in\ker\ad_{e_0}$ then $\mathbf{L}_0=\mathbf{W}\rtimes\F e_1=\mathbf{W}\oplus\F e_1$ is Abelian. Otherwise $\mathbf{L}_0=\mathbf{W}\rtimes\F e_1$ is almost Abelian.

Now let $\mathbf{I}=\mathbf{L}_0\subset\mathbf{L}$ be a Lie subalgebra as above. In case $\mathbf{I}=\mathbf{W}\subset\mathbf{V}$ then $\mathbf{I}$ is an ideal when
$$
[\mathbf{W},\mathbf{L}]=\ad_{e_0}\mathbf{W}\subset\mathbf{W}.
$$
Otherwise if $\mathbf{I}=\mathbf{W}\rtimes\F e_1$ then it is an ideal when
$$
[\mathbf{I},\mathbf{L}]=[\mathbf{W}\rtimes\F e_1,\mathbf{V}\rtimes\F e_0]=\ad_{e_0}\mathbf{V}\subset\mathbf{W}.
$$
This completes the proof.
\end{proof}

\begin{remark} Note that the second kind of Abelian ideals in the proposition above is possible only if $\mathbf{L}$ is 2-step nilpotent.
\end{remark}

In later sections we will observe that the direct sum $\mathbf{L}\oplus\mathbf{W}$ of an almost Abelian Lie algebra $\mathbf{L}$ with an Abelian Lie algebra $\mathbf{W}$ is not qualitatively different from $\mathbf{L}$, hence it makes sense to separate out Abelian direct factor subalgebras which are called central extensions.

\begin{lemma}\label{aADirectDecLemma} If an almost Abelian Lie algebra is decomposable,
$$
\mathbf{L}=\mathbf{L}_1\oplus\mathbf{L}_2,
$$
then $\mathbf{L}_1$ is almost Abelian and $\mathbf{L}_2$ is Abelian or vice versa.
\end{lemma}
\begin{proof} Let $\operatorname{P}_1$ and $\operatorname{P}_2$ be projectors onto the subspaces $\mathbf{L}_1$ and $\mathbf{L}_2$, $\operatorname{P}_1+\operatorname{P}_2=\id$. Let further $\mathbf{L}=\mathbf{V}\rtimes\F e_0$ as usual. Then $\operatorname{P}_1\mathbf{V}\subset\mathbf{L}_1$, $\operatorname{P}_2\mathbf{V}\subset\mathbf{L}_2$ and
$$
\mathbf{V}=(\operatorname{P}_1+\operatorname{P}_2)\mathbf{V}\subset\operatorname{P}_1\mathbf{V}+\operatorname{P}_2\mathbf{V}=\operatorname{P}_1\mathbf{V}\oplus\operatorname{P}_2\mathbf{V}.
$$
It follows that $\codim_\F(\operatorname{P}_1\mathbf{V}+\operatorname{P}_2\mathbf{V})\le1$ and therefore, without loss of generality, $\operatorname{P}_2\mathbf{V}=\mathbf{L}_2$ and $\codim_\F\operatorname{P}_1\mathbf{V}\le1$ in $\mathbf{L}_1$. Now $\mathbf{L}=\mathbf{L}_1\oplus\mathbf{L}_2$ implies that $\operatorname{P}_1:\mathbf{L}\to\mathbf{L}_1$ and $\operatorname{P}_2:\mathbf{L}\to\mathbf{L}_2$ are Lie algebra homomorphisms, thus, both $\operatorname{P}_2\mathbf{V}=\mathbf{L}_2$ and $\operatorname{P}_1\mathbf{V}$ are Abelian subalgebras whence the assertion follows.
\end{proof}

\begin{remark}\label{LDecompRemark} If an almost Abelian Lie algebra decomposes as $\mathbf{L}=\mathbf{L}_0\oplus_\F\mathbf{W}$ for two subspaces $\mathbf{L}_0,\mathbf{W}\subset\mathbf{L}$ such that
$[\mathbf{L},\mathbf{L}]\subset\mathbf{L}_0$ and $\mathbf{W}\subset\mathcal{Z}(\mathbf{L})$ then $\mathbf{L}=\mathbf{L}_0\oplus\mathbf{W}$. Thus, if $\mathbf{L}$ is indecomposable then $\mathcal{Z}(\mathbf{L})\subset[\mathbf{L},\mathbf{L}]$.
\end{remark}

\begin{proposition}\label{DecompProp} Every almost Abelian Lie algebra $\mathbf{L}$ can be written as
$$
\mathbf{L}=\mathbf{L}_0\oplus\mathbf{W},
$$
where $\mathbf{L}_0$ is an indecomposable almost Abelian Lie subalgebra and $\mathbf{W}$ is a possibly trivial Abelian Lie subalgebra. If
$$
\mathbf{L}=\mathbf{L}_0'\oplus\mathbf{W}'
$$
is another such decomposition then $\mathbf{L}_0\simeq\mathbf{L}_0'$ and $\mathbf{W}\simeq\mathbf{W}'$.
\end{proposition}
\begin{proof}
Choose as $\mathbf{W}$ any complement of $[\mathbf{L},\mathbf{L}]\cap\mathcal{Z}(\mathbf{L})$ in $\mathcal{Z}(\mathbf{L})$, i.e.,
\begin{equation}
\mathcal{Z}(\mathbf{L})=\left([\mathbf{L},\mathbf{L}]\cap\mathcal{Z}(\mathbf{L})\right)\oplus_\F\mathbf{W}.\label{kerade0Decomp}
\end{equation}
Then by Remark \ref{LDecompRemark} we get $\mathbf{L}=\mathbf{L}_0\oplus\mathbf{W}$ for some Lie subalgebra $\mathbf{L}_0\subset\mathbf{L}$, and by Lemma \ref{aADirectDecLemma} we have that $\mathbf{L}_0$ is almost Abelian. If $\mathbf{L}_0$ is decomposable then $\mathbf{L}_0=\mathbf{L}_1\oplus\mathbf{W}_1$ with a non-trivial Abelian subalgebra $\mathbf{W}_1\neq 0$. It follows that $\mathbf{W}\oplus\mathbf{W}_1\subset\mathcal{Z}(\mathbf{L})$ but from (\ref{kerade0Decomp}) we find that $(\mathbf{W}\oplus\mathbf{W}_1)\cap[\mathbf{L},\mathbf{L}]\neq0$ which is a contradiction.

Suppose that $\mathbf{L}=\mathbf{L}_0'\oplus\mathbf{W}'$ is another decomposition. From Remark \ref{LDecompRemark} we find that
$$
\mathcal{Z}(\mathbf{L})=\mathcal{Z}(\mathbf{L}_0)\oplus\mathbf{W}=(\mathcal{Z}(\mathbf{L}_0)\cap[\mathbf{L}_0,\mathbf{L}_0])\oplus\mathbf{W}=(\mathcal{Z}(\mathbf{L})\cap[\mathbf{L},\mathbf{L}])\oplus\mathbf{W}=(\mathcal{Z}(\mathbf{L})\cap[\mathbf{L},\mathbf{L}])\oplus\mathbf{W}',
$$
hence $\mathbf{W}\simeq\mathbf{W}'$. But then
$$
\mathbf{L}_0\oplus\mathbf{W}=\mathbf{L}_0'\oplus\mathbf{W}'\simeq\mathbf{L}_0'\oplus\mathbf{W},
$$
whence $\mathbf{L}_0\simeq\mathbf{L}_0'$ by dividing out $\mathbf{W}$ on both sides.
\end{proof}
We conclude this section by introducing some of the most prominent almost Abelian Lie algebras as examples.
\begin{example}\label{HeisExamp} Denote by
$$
\mathbf{H}_\F=\left\{\begin{pmatrix}
0 & 0 & 0\\
t & 0 & p\\
q & 0 & 0
\end{pmatrix}\vline\quad(p,t,q)\in\F^3\right\}
$$
the Heisenberg algebra which is special in many respects. This corresponds to $\mathbf{H}_\F=\F^2\rtimes\F e_0$ with
$$
\ad_{e_0}=\begin{pmatrix}
0 & 1\\
0 & 0
\end{pmatrix}.
$$
\end{example}
\begin{example} Denote by
$$
\mathbf{ax\!+\!b}_\F=\left\{\begin{pmatrix}
0 & 0\\
b & a
\end{pmatrix}\vline\quad(a,b)\in\F^2\right\}
$$
the Lie algebra of generators of affine transformations on $\F$. This is up to isomorphism the only non-Abelian 2-dimensional Lie algebra. This corresponds to $\mathbf{ax\!+\!b}_\F=\F\rtimes\F e_0$ with $\ad_{e_0}=\id$.
\end{example}


\section{Homomorphisms and classification\label{HomClass}}

This section is devoted to homomorphisms and identification or classification of almost Abelian Lie algebras. We start by showing that for an almost Abelian Lie algebra which is not a central extension of the Heisenberg algebra (i.e.,  $\mathbf{L}\not\simeq\mathbf{H}_\F\oplus\mathbf{W}$) the form $\mathbf{L}=\mathbf{V}\rtimes\F e_0$ is unique.

\begin{proposition}\label{2cod1AbIdProp} If an almost Abelian Lie algebra $\mathbf{L}$ has more than one codimension 1 Abelian ideals then $\mathbf{L}\simeq\mathbf{H}_\F\oplus\mathbf{W}$ with $\mathbf{W}$ Abelian.
\end{proposition}
\begin{proof} Let $\mathbf{V}\neq\mathbf{V}'$ be two distinct codimension 1 Abelian ideals so that $\mathbf{L}=\mathbf{V}\rtimes\F e_0=\mathbf{V}'\rtimes\F e_0'$. Denote $\mathbf{V}''\doteq\mathbf{V}\cap\mathbf{V}'$. Then there is an element $v_1\in\mathbf{V}$ ($v_1'\in\mathbf{V}'$) such that $v_1\not\in\mathbf{V}'$ ($v_1'\not\in\mathbf{V}$) and $\mathbf{V}=\F v_1\oplus_\F\mathbf{V}''$ ($\mathbf{V}'=\F v_1'\oplus_\F\mathbf{V}''$). That $v_1\not\in\mathbf{V}'$ implies that $v_1=\lambda e_0'+w_1$ for some $\lambda\in\F^*$ and $w_1\in\mathbf{V}'$. It follows that
$$
[e_0',\mathbf{V}'']=\frac1\lambda[v_1-w_1,\mathbf{V}'']=0,
$$
i.e., $\mathbf{V}''\subset\mathcal{Z}(\mathbf{L})$. Because $\mathbf{L}$ is non-Abelian we know that $v_1\notin\mathcal{Z}(\mathbf{L})$, thus $\mathbf{V}''=\mathcal{Z}(\mathbf{L})$. Denoting $v_2\doteq[e_0,v_1]\neq0$ we find that
$$
[\mathbf{L},\mathbf{L}]=[e_0,\mathbf{V}]=[e_0,\F v_1\oplus_\F\mathbf{V}'']=\F v_2.
$$
Similarly, if we denote $v_2'\doteq[e_0',v_1']\neq0$ then $[\mathbf{L},\mathbf{L}]=\F v_2'$. It follows that $v_2'=\mu v_2$ for some $\mu\in\F^*$. Now $v_2\in\mathbf{V}$ and $v_2=\mu^{-1}v_2'\in\mathbf{V}'$ therefore $v_2\in\mathbf{V}\cap\mathbf{V}'=\mathbf{V}''$. Let $\mathbf{V}''=\F v_2\oplus\mathbf{W}$. Then by Remark \ref{LDecompRemark} we get
$$
\mathbf{L}=\mathbf{L}_0\oplus\mathbf{W},\quad\mathbf{L}_0\doteq(\F v_1\oplus_\F\F v_2)\rtimes\F e_0.
$$
Now it is straightforward to verify that
$$
\mathbf{L}_0\ni pe_0+qv_1+tv_2\longleftrightarrow\begin{pmatrix}
0 & 0 & 0\\
t & 0 & p\\
q & 0 & 0
\end{pmatrix}\in\mathbf{H}_\F,\quad(p,q,t)\in\F^3
$$
is a Lie algebra isomorphism.
\end{proof}

\begin{remark} It can be observed from Proposition \ref{SubalgIdealProp} and Remark \ref{ZRemark} that $\mathbf{L}$ has a codimension 1 Abelian ideal other than $\mathbf{V}$ if and only if $\dim\mathbf{V}/\ker\ad_{e_0}=\dim\mathbf{V}/\mathcal{Z}(\mathbf{L})=1$ and $\ad_{e_0}\mathbf{V}=[\mathbf{L},\mathbf{L}]\subset\mathcal{Z}(\mathbf{L})$.
\end{remark}
In principle, Proposition \ref{2cod1AbIdProp} could be given a simpler proof based on the above remark. In fact, a stronger statement can be obtained easily in the same spirit.

\begin{remark} If an almost Abelian Lie algebra $\mathbf{L}$ has more than one codimension 1 Abelian subalgebras then $\mathbf{L}\simeq\mathbf{H}_\F\oplus\mathbf{W}$ or $\mathbf{L}\simeq\mathbf{ax\!+\!b}_\F\oplus\mathbf{W}$ with $\mathbf{W}$ Abelian.
\end{remark}
Indeed, it can be observed from Proposition \ref{SubalgIdealProp} and Remark \ref{ZRemark} that $\mathbf{L}=\mathbf{L}_0\oplus\mathbf{W}$ has a codimension 1 Abelian subalgebra other than $\mathbf{V}$ if and only if $\dim\mathbf{V}/\mathcal{Z}(\mathbf{L})=1$, in which case $\dim[\mathbf{L}_0,\mathbf{L}_0]=\dim\ad_{e_0}\mathbf{V}=1$, and then we use Lemma \ref{1dimadVLemma} from below.

We continue the study of almost Abelian Lie algebras by considering structure preserving maps between them. Our first subject is a Lie algebra homomorphism $\phi:\mathbf{L}\to\mathbf{L}'$ from an almost Abelian Lie algebra $\mathbf{L}$ to another Lie algebra $\mathbf{L}'$. The image $\phi(\mathbf{L})\subset\mathbf{L}'$ is a Lie subalgebra, and the rest of $\mathbf{L}'$ is completely irrelevant from the point of view of the homomorphism $\phi$. Therefore in studying $\phi$ it is no loss of generality to assume that it is onto, i.e., $\phi(\mathbf{L})=\mathbf{L}'$. Denote by $\Hom(\mathbf{L},\mathbf{L}')$ the $\F$-vector space of all Lie algebra homomorphisms $\phi:\mathbf{L}\to\mathbf{L}'$ and by $\Hom(\mathbf{L}\onto\mathbf{L}')\subset\Hom(\mathbf{L},\mathbf{L}')$ the subset of surjective homomorphisms.

\begin{remark}\label{HomImRemark} The image of a homomorphism $\phi:\mathbf{L}\to\mathbf{L}'$ from an almost Abelian Lie algebra $\mathbf{L}=\mathbf{V}\rtimes\F e_0$ is either Abelian or almost Abelian. Indeed, $\phi\mathbf{V}\subset\phi\mathbf{L}$ is obviously an Abelian ideal with $\codim\phi\mathbf{V}\le1$.
\end{remark}
Consider a surjective homomorphism $\phi:\mathbf{L}\to\mathbf{L}'$. According to the Fundamental Theorem on Lie algebra homomorphisms (an explanation can be found in \cite{Bou98}) there exists a unique Lie algebra isomorphism $\psi:\mathbf{L}/\ker\phi\to\mathbf{L}'$ such that $\phi=\psi\circ\mathfrak{q}$ where $\mathfrak{q}:\mathbf{L}\to\mathbf{L}/\ker\phi$ is the canonical quotient homomorphism. Conversely, every isomorphism $\psi:\mathbf{L}/\ker\phi\to\mathbf{L}'$ composed with $\mathfrak{q}$ gives a homomorphism $\psi\circ\mathfrak{q}:\mathbf{L}\to\mathbf{L}'$ with the same kernel $\ker\phi$. Moreover, given a fixed isomorphism $\psi_0:\mathbf{L}/\ker\phi\to\mathbf{L}'$ every other isomorphism $\psi:\mathbf{L}/\ker\phi\to\mathbf{L}'$ can be written as $\psi=\psi_0\circ\Psi$ where $\Psi\in\Aut(\mathbf{L}/\ker\phi)$ is an automorphism. And every automorphism $\Psi$ gives rise to a new isomorphism $\psi=\psi_0\circ\Psi$. We conclude in the following remark.

\begin{remark}\label{HomDescRemark} There is a bijective correspondence between surjective homomorphisms $\phi:\mathbf{L}\to\mathbf{L}'$ from a fixed Lie algebra $\mathbf{L}$ and pairs $(\mathbf{I},\Psi)$ where $\mathbf{I}=\ker\phi\subset\mathbf{L}$ is an ideal and $\Psi\in\Aut(\mathbf{L}/\mathbf{I})$.
\end{remark}

Let $\phi:\mathbf{L}\to\mathbf{L}'$ be a Lie algebra homomorphism from an almost Abelian Lie algebra $\mathbf{L}$ into an Abelian or almost Abelian Lie algebra $\mathbf{L}'$. According to Proposition \ref{DecompProp} let $\mathbf{L}=\mathbf{L}_0\oplus\mathbf{W}$ with $\mathbf{L}_0$ an indecomposable almost Abelian and $\mathbf{W}$ a possibly trivial Abelian Lie algebra. In a similar fashion let $\mathbf{L}'=\mathbf{L}_0'\oplus\mathbf{W}'$ where $\mathbf{L}_0'$ is either a trivial or an indecomposable almost Abelian and $\mathbf{W}'$ is a possibly trivial Abelian Lie algebra. Then $\phi$ can be broken down into the block form
\begin{equation}
\begin{pmatrix}
\mathbf{L}_0'\\
\mathbf{W}'
\end{pmatrix}=\begin{pmatrix}
\phi_{00} & \phi_{01}\\
\phi_{10} & \phi_{11}
\end{pmatrix}\begin{pmatrix}
\mathbf{L}_0\\
\mathbf{W}
\end{pmatrix}\label{phiBlockForm}
\end{equation}

\begin{proposition}\label{SurjHomProp} The following describes the block form (\ref{phiBlockForm}) of a surjective homomorphism $\phi\in\Hom(\mathbf{L}\onto\mathbf{L}')$ from an almost Abelian Lie algebras $\mathbf{L}$ onto an Abelian or almost Abelian Lie algebra $\mathbf{L}'$:
\begin{itemize}
\item $\phi_{00}\in\Hom(\mathbf{L}_0\onto\mathbf{L}_0')$

\item $\phi_{01}\in\Hom(\mathbf{W},\mathcal{Z}(\mathbf{L}_0'))$

\item $\phi_{10}\in\Hom(\mathbf{L}_0,\mathbf{W}')$, i.e., $[\mathbf{L}_0,\mathbf{L}_0]\subset\ker\phi_{10}$

\item $\phi_{11}\in\Hom(\mathbf{W},\mathbf{W}')$ such that $\phi_{10}\mathbf{L}_0+\phi_{11}\mathbf{W}=\mathbf{W}'$
\end{itemize}
\end{proposition}
\begin{proof} By surjectivity of $\phi$ we have $\mathbf{L}_0'=\phi_{00}\mathbf{L}_0+\phi_{01}\mathbf{W}$. Then using $[\mathbf{L}_0,\mathbf{W}]=0$ and the homomorphism condition for $\phi$ we find that $\phi_{01}\mathbf{W}\subset\mathcal{Z}(\mathbf{L}_0')$. But because $\mathbf{L}_0'$ is indecomposable by Remark \ref{LDecompRemark} we have
$$
\mathcal{Z}(\mathbf{L}_0')\subset[\mathbf{L}_0',\mathbf{L}_0']=[\phi_{00}\mathbf{L}_0,\phi_{00}\mathbf{L}_0]\subset\phi_{00}\mathbf{L}_0,
$$
therefore $\mathbf{L}_0'=\phi_{00}\mathbf{L}_0+\phi_{01}\mathbf{W}\subset\phi_{00}\mathbf{L}_0$ proving the surjectivity of $\phi_{00}$. The rest is an easy consequence of Remark \ref{HomDescRemark} and the preceding discussion.
\end{proof}
Consider now a homomorphism onto itself, i.e., $\phi\in\Hom(\mathbf{L}\onto\mathbf{L})$.
\begin{lemma}\label{phiInjphi00InjLemma} If a surjective homomorphism $\phi\in\Hom(\mathbf{L}\onto\mathbf{L})$ is injective, $\ker\phi=0$, then so is its component $\phi_{00}\in\Hom(\mathbf{L}_0\onto\mathbf{L}_0)$, i.e., $\ker\phi_{00}=0$.
\end{lemma}
\begin{proof} If $\ker\phi=0$ then the only solution of the equation
$$
\phi X=\begin{pmatrix}
\phi_{00} & \phi_{01}\\
\phi_{10} & \phi_{11}
\end{pmatrix}\begin{pmatrix}
X\\
0
\end{pmatrix}=\begin{pmatrix}
\phi_{00}X\\
\phi_{10}X
\end{pmatrix}=0,\quad X\in\mathbf{L}_0
$$
is $X=0$, which means $\ker\phi_{00}\cap\ker\phi_{10}=0$. Because $\ker\phi_{00}\subset\mathbf{L}_0$ is an ideal this implies
\begin{equation}
[\mathbf{L}_0,\ker\phi_{00}]\subset\ker\phi_{00}\cap[\mathbf{L}_0,\mathbf{L}_0]\subset\ker\phi_{00}\cap\ker\phi_{10}=0,\label{kerphi00}
\end{equation}
where we used $[\mathbf{L}_0,\mathbf{L}_0]\subset\ker\phi_{10}$ from Proposition \ref{SurjHomProp}. Thus , $\ker\phi_{00}\subset\mathcal{Z}(\mathbf{L}_0)$. But $\mathbf{L}_0$ is indecomposable, hence
$$
\ker\phi_{00}\subset\mathcal{Z}(\mathbf{L}_0)\subset[\mathbf{L}_0,\mathbf{L}_0],
$$
and thereby using (\ref{kerphi00}) we get
$$
\ker\phi_{00}=\ker\phi_{00}\cap[\mathbf{L}_0,\mathbf{L}_0]=0
$$
which completes the proof.
\end{proof}
\begin{lemma}\label{phi00InjkerphiZLemma} For a surjective homomorphism $\phi\in\Hom(\mathbf{L}\onto\mathbf{L})$, if the component $\phi_{00}\in\Hom(\mathbf{L}_0\onto\mathbf{L}_0)$ is injective, $\ker\phi_{00}=0$, then $\ker\phi\subset\mathcal{Z}(\mathbf{L}_0)\oplus\mathbf{W}$.
\end{lemma}
\begin{proof} Let $X+Y\in\ker\phi\subset\mathbf{L}=\mathbf{L}_0\oplus\mathbf{W}$, that is,
$$
\phi(X+Y)=\begin{pmatrix}
\phi_{00} & \phi_{01}\\
\phi_{10} & \phi_{11}
\end{pmatrix}\begin{pmatrix}
X\\
Y
\end{pmatrix}=\begin{pmatrix}
\phi_{00}X+\phi_{01}Y\\
\phi_{10}X+\phi_{11}Y
\end{pmatrix}=0,\quad X\in\mathbf{L}_0,\quad Y\in\mathbf{W}.
$$
By Proposition \ref{SurjHomProp} we have $\phi_{00}X=-\phi_{01}Y\in\mathcal{Z}(\mathbf{L}_0)$. This means
$$
0=[\phi_{00}X,\mathbf{L}_0]=[\phi_{00}X,\phi_{00}\mathbf{L}_0]=\phi_{00}[X,\mathbf{L}_0]
$$
whence
$$
[X,\mathbf{L}_0]\in\ker\phi_{00}=0,
$$
that is, $X\in\mathcal{Z}(\mathbf{L}_0)$.
\end{proof}
\begin{lemma}\label{kerphiZphiInjphi00phi11InjLemma} Let $\phi\in\Hom(\mathbf{L}\onto\mathbf{L})$ with $\ker\phi\subset\mathcal{Z}(\mathbf{L}_0)\oplus\mathbf{W}$. Then $\ker\phi=0$ if and only if $\ker\phi_{00}\cap\mathcal{Z}(\mathbf{L}_0)=0$ and $\ker\phi_{11}=0$.
\end{lemma}
\begin{proof} As $\ker\phi\subset\mathcal{Z}(\mathbf{L}_0)\oplus\mathbf{W}$ we have that $\phi$ is injective if and only if its restriction to $\mathcal{Z}(\mathbf{L}_0)\oplus\mathbf{W}$ is. Because $\mathbf{L}_0$ is indecomposable and thus, by Remark \ref{LDecompRemark}
$$
\mathcal{Z}(\mathbf{L}_0)\subset[\mathbf{L}_0,\mathbf{L}_0]\subset\ker\phi_{10},
$$
the restriction of $\phi$ to $\mathcal{Z}(\mathbf{L}_0)\oplus\mathbf{W}$ has the form
$$
\begin{pmatrix}
\phi_{00} & \phi_{01}\\
\phi_{10} & \phi_{11}
\end{pmatrix}\begin{pmatrix}
\mathcal{Z}(\mathbf{L}_0)\\
\mathbf{W}
\end{pmatrix}=\begin{pmatrix}
\phi_{00} & \phi_{01}\\
0 & \phi_{11}
\end{pmatrix}\begin{pmatrix}
\mathcal{Z}(\mathbf{L}_0)\\
\mathbf{W}
\end{pmatrix}=\begin{pmatrix}
\id & 0\\
0 & \phi_{11}
\end{pmatrix}\begin{pmatrix}
\id & \phi_{01}\\
0 & \id
\end{pmatrix}\begin{pmatrix}
\phi_{00} & 0\\
0 & \id
\end{pmatrix}\begin{pmatrix}
\mathcal{Z}(\mathbf{L}_0)\\
\mathbf{W}
\end{pmatrix}.
$$
From here it is clear that $\ker\phi=0$ if and only if $\ker\phi_{00}\cap\mathcal{Z}(\mathbf{L}_0)=0$ and $\ker\phi_{11}=0$.
\end{proof}
We are finally ready to describe the automorphism group $\Aut(\mathbf{L})$ of an almost Abelian Lie algebra $\mathbf{L}$.

\begin{proposition} The automorphism group of an almost Abelian Lie algebra $\mathbf{L}=\mathbf{L}_0\oplus\mathbf{W}$ with $\mathbf{L}_0$ indecomposable has the following form,
$$
\Aut(\mathbf{L})=\left\{\begin{pmatrix}
\phi_{00} & \phi_{01}\\
\phi_{10} & \phi_{11}
\end{pmatrix}\in\Hom(\mathbf{L}\onto\mathbf{L})\vline\quad\phi_{00}\in\Aut(\mathbf{L}_0),\quad \ker\phi_{11}=0\right\}.
$$
\end{proposition}
\begin{proof} Let $\phi\in\Aut(\mathbf{L})$, i.e., $\phi\in\Hom(\mathbf{L}\onto\mathbf{L})$ and $\ker\phi=0$. Then by Lemma \ref{phiInjphi00InjLemma} $\ker\phi_{00}=0$ (i.e., $\phi_00\in\Aut(\mathbf{L}_0)$) and by Lemma \ref{kerphiZphiInjphi00phi11InjLemma} $\ker\phi_{11}=0$. Conversely, let $\phi\in\Hom(\mathbf{L}\onto\mathbf{L})$ such that $\ker\phi_{00}=\ker\phi_{11}=0$. Then by Lemma \ref{phi00InjkerphiZLemma} we have $\ker\phi\subset\mathcal{Z}(\mathbf{L}_0)\oplus\mathbf{W}$, and then by Lemma \ref{kerphiZphiInjphi00phi11InjLemma} we establish that $\ker\phi=0$.
\end{proof}
Homomorphisms between Abelian Lie algebras are simply linear operators. In Proposition \ref{SubalgIdealProp} we have classified all subalgebras and ideals of almost Abelian Lie algebras. By the above proposition we only need to study $\Hom(\mathbf{L}_0\onto\mathbf{L}_0')$ in order to have a full understanding of $\Hom(\mathbf{L},\mathbf{L}')$.

\begin{remark}\label{phiLsimL1Remark} If $\phi\in\Hom(\mathbf{L}\onto\mathbf{L}')$ and $\mathbf{L}=\mathbf{V}\rtimes\F e_0$ then $\mathbf{L}'=\phi\mathbf{V}\rtimes\F\phi e_0$ and $\phi\circ\ad_{e_0}=\ad_{\phi e_0}\circ\phi$.
\end{remark}

By Remark \ref{HomDescRemark} all we need to do is to study $\Aut(\mathbf{L})$ for indecomposable almost Abelian Lie algebras $\mathbf{L}$ (for Abelian Lie algebras $\Aut(\mathbf{L})$ is the group of invertible operators). For an almost Abelian Lie algebra $\mathbf{L}=\mathbf{V}\rtimes\F e_0$ the action of an endomorphism $\phi\in\End_\F(\mathbf{L})$ on an element $X=(v,t)\in\mathbf{V}\oplus_\F\F$ can be written in the block form
\begin{equation}
\phi X=\begin{pmatrix}
\Delta && \gamma\\
\beta^\top && \alpha
\end{pmatrix}\begin{pmatrix}
v\\
t
\end{pmatrix},\quad\alpha\in\F,\quad\beta^\top\in\mathbf{V}^*,\quad\gamma\in\mathbf{V},\quad\Delta\in\End_\F(\mathbf{V}).\label{EndBlockForm}
\end{equation}
By $\mathbf{V}^*$ (unlike $\F^*$ for a field) we denote the space of linear functionals on $\mathbf{V}$, i.e., the dual space. The automorphism group of the Heisenberg algebra $\Aut(\mathbf{H}_\F)$ is a classical subject (e.g., \cite{Fol89}, \cite{Tha98}) and we include a proposition merely for completeness.

\begin{proposition}\label{HeisAutProp} The automorphism group of the Heisenberg algebra in the coordinates $(p,t,q)$ of Example \ref{HeisExamp} is
$$
\Aut(\mathbf{H}_\F)=\left\{\begin{pmatrix}
\alpha\Delta_{22}-\beta_2\gamma_2 & \Delta_{12} & \gamma_1\\
0 & \Delta_{22} & \gamma_2\\
0 & \beta_2 & \alpha
\end{pmatrix}\vline\quad\alpha,\beta_2,\gamma_1,\gamma_2,\Delta_{12},\Delta_{22}\in\F,\quad\alpha\Delta_{22}-\beta_2\gamma_2\neq0\right\}.
$$
\end{proposition}
\begin{proof} Simply plug in the block form (\ref{EndBlockForm}) into the homomorphism condition $\phi[X,Y]=[\phi X,\phi Y]$ for all $X,Y\in\mathbf{H}_\F$, then require invertibility of $\phi$.
\end{proof}
Automorphism groups of all other indecomposable almost Abelian Lie algebras are more restricted in view of Proposition \ref{2cod1AbIdProp}.

\begin{proposition} The automorphism group of an indecomposable almost Abelian Lie algebra $\mathbf{L}=\mathbf{V}\rtimes\F e_0$ other than $\mathbf{H}_\F$ is
$$
\Aut(\mathbf{L})=\left\{\begin{pmatrix}
\Delta & \gamma\\
0 & \alpha
\end{pmatrix}\vline\quad\alpha\in\F^*,\quad\gamma\in\mathbf{V},\quad\Delta\in\Aut(\mathbf{V}),\quad\Delta\ad_{e_0}-\alpha\ad_{e_0}\Delta=0\right\}.
$$
\end{proposition}
\begin{proof} Let $\phi:\mathbf{V}\to\mathbf{V}$ be an automorphism in the block form (\ref{EndBlockForm}). Then $\phi v=\beta^\top v e_0+\Delta v$ for $\forall v\in\mathbf{V}$. Because $\phi$ is bijective we know that $\phi\mathbf{V}\subset\mathbf{L}$ is a codimension 1 Abelian ideal, and as $\mathcal{L}$ is not isomorphic to the Heisenberg algebra we conclude by Proposition \ref{2cod1AbIdProp} that $\phi\mathbf{V}=\mathbf{V}$. Thus $\beta^\top=0$. Write
$$
\phi=\begin{pmatrix}
\Delta & \gamma\\
0 & \alpha
\end{pmatrix}=\begin{pmatrix}
\id & 0\\
0 & \alpha
\end{pmatrix}\begin{pmatrix}
\id & \gamma\\
0 & 1
\end{pmatrix}\begin{pmatrix}
\Delta & 0\\
0 & 1
\end{pmatrix}.
$$
$\phi$ is bijective if and only if all three factors are, which is equivalent to $\alpha\neq 0$ and $\Delta$ being bijective. Finally the homomorphism condition for $\phi$ reads
$$
\phi\ad_{e_0}v=\Delta\ad_{e_0}v=\phi[e_0,v]=[\phi e_0,\phi v]=[\alpha e_0+\gamma,\Delta v]=\alpha\ad_{e_0}\Delta v,\quad\forall v\in\mathbf{V},
$$
precisely as in the statement.
\end{proof}

\begin{corollary}\label{AutV1V2Corr} For any two codimension 1 Abelian ideals $\mathbf{V},\mathbf{V}'\subset\mathbf{L}$ of an almost Abelian Lie algebra $\mathbf{L}$ there exists an automorphism $\phi\in\Aut(\mathbf{L})$ such that $\phi\mathbf{V}=\mathbf{V}'$.
\end{corollary}
\begin{proof} By Proposition \ref{2cod1AbIdProp} either $\mathbf{V}=\mathbf{V'}$ and we take $\phi=\id$ or $\mathbf{L}=\mathbf{H}_\F\oplus\mathbf{W}$ which we will assume. If $\mathbf{L}=\mathbf{V}\rtimes\F e_0$ then we write $\mathbf{H}_\F=(\F v_1\oplus_\F\F v_2)\rtimes\F e_0$ with $[e_0,v_1]=v_2$ and $[e_0,v_2]=0$. Thus $\mathbf{V}=\F v_1\oplus\F v_2\oplus\mathbf{W}$. By Proposition \ref{SubalgIdealProp} a codimension 1 Abelian ideal other than $\mathbf{V}$ must have the form $\mathbf{V'}=\F e_1\oplus\F v_2\oplus\mathbf{W}$ with $e_1=e_0+\lambda v_1$ for some $\lambda\in\F$. In the basis $(e_0,v_2,v_1)$ of Proposition \ref{HeisAutProp} we have
$$
\F v_1\oplus\F v_2=\left\{\begin{pmatrix}
0\\
\mu\\
\nu
\end{pmatrix},\quad(\mu,\nu)\in\F^2\right\},\quad\F e_1\oplus\F v_2=\left\{\begin{pmatrix}
\eta\\
\rho\\
\lambda\eta
\end{pmatrix},\quad(\eta,\rho)\in\F^2\right\}.
$$
Now it is easy to check that the automorphism $\phi=\psi\oplus\id$ with
$$
\psi=\begin{pmatrix}
1 & 0 & 1\\
0 & 1 & 0\\
\lambda-1 & 0 & \lambda
\end{pmatrix}
$$
will do the job.
\end{proof}

Our final task will be the classification of almost Abelian Lie algebras into isomorphism classes. Clearly, two Lie algebras $\mathbf{L}=\mathbf{L}_0\oplus\mathbf{W}$ and $\mathbf{L}'=\mathbf{L_0}'\oplus\mathbf{W}'$ are isomorphic if and only if $\mathbf{L}_0\simeq\mathbf{L}_0'$ and $\mathbf{W}\simeq\mathbf{W}'$. It turns out that isomorphism classes of indecomposable almost Abelian Lie algebras $\mathbf{L}=\mathbf{V}\rtimes\F e_0$ correspond to similarity classes of the operators $\ad_{e_0}$ up to rescaling. Recall the description of $\mathbf{L}=\mathbf{V}\rtimes\F e_0$ in terms of the pair $(\mathbf{V},\ad_{e_0})$ of Remark \ref{VadPairRemark}.

\begin{definition} Two pairs $(\mathbf{V}_1,\mbox{\normalfont T}_1)$ and $(\mathbf{V}_2,\mbox{\normalfont T}_2)$ with $\mbox{\normalfont T}_1\in\End_\F(\mathbf{V}_1)$ and $\mbox{\normalfont T}_2\in\End_\F(\mathbf{V}_2)$ are called similar, $(\mathbf{V}_1,\mbox{\normalfont T}_1)\sim(\mathbf{V}_2,\mbox{\normalfont T}_2)$, if there exists an invertible map $\phi:\mathbf{V}_1\to\mathbf{V}_2$ such that $\phi\mbox{\normalfont T}_1=\mbox{\normalfont T}_2\phi$.
\end{definition}
Similarity is an equivalence relation in the class of pairs $(\mathbf{V},\mbox{T})$. For every such pair with $\mbox{T}\neq0$ we can construct the almost Abelian Lie algebra $\mathbf{L}=\mathbf{V}\rtimes\F e_0$ where $\ad_{e_0}=\mbox{T}$. Conversely, every almost Abelian Lie algebra $\mathbf{L}=\mathbf{V}\rtimes\F e_0$ gives a pair $(\mathbf{V},\ad_{e_0})$ with $\ad_{e_0}\neq0$.

\begin{proposition} Two almost Abelian Lie algebras $\mathbf{L}=\mathbf{V}\rtimes\F e_0$ and $\mathbf{L}'=\mathbf{V}'\rtimes\F e_0'$ are isomorphic if and only if $(\mathbf{V},\ad_{e_0})\sim(\mathbf{V}',\lambda\ad_{e_0'})$ for some $\lambda\in\F^*$.
\end{proposition}
\begin{proof} Assume first that $(\mathbf{V},\ad_{e_0})\sim(\mathbf{V}',\lambda\ad_{e_0'})$, i.e., there is an invertible $\phi:\mathbf{V}\to\mathbf{V}'$ such that $\phi\ad_{e_0}=\lambda\ad_{e_0'}\phi$. Define the linear invertible map $\Phi:\mathbf{L}\to\mathbf{L}'$ by setting $\Phi e_0=\lambda e_0'$ and $\Phi v=\phi v$ for $v\in\mathbf{V}$. It remains to note that $\Phi$ is a Lie algebra isomorphism, because
$$
\Phi[e_0,v]=\Phi\ad_{e_0}v=\phi\ad_{e_0}v=\lambda\ad_{e_0'}\phi v=[\lambda e_0',\phi v]=[\Phi e_0,\Phi v],\quad\forall v\in\mathbf{V}.
$$
Conversely, let $\Phi:\mathbf{L}\to\mathbf{L}'$ be a Lie algebra isomorphism. Then $\Phi\mathbf{V}\subset\mathbf{L}'$ is a codimension 1 Abelian ideal, and by Corollary \ref{AutV1V2Corr} we have an automorphism $\Psi\in\Aut(\mathbf{L}')$ such that $\Psi\Phi\mathbf{V}=\mathbf{V}'$. Now $\Psi\Phi:\mathbf{L}\to\mathbf{L}'$ is another isomorphism. By Remark \ref{phiLsimL1Remark} we have
$$
\mathbf{L}'=\mathbf{V}'\rtimes\F e_0'=\Psi\Phi\mathbf{V}\rtimes\F\Psi\Phi e_0=\mathbf{V}'\rtimes\F\Psi\Phi e_0,\quad\Psi\Phi\ad_{e_0}=\ad_{\Psi\Phi e_0}\Psi\Phi.
$$
Because $e_0'\notin\Psi\Phi\mathbf{V}=\mathbf{V}'$ we have that $e_0'=\lambda\Psi\Phi e_0+v_0'$ for some $\lambda\in\F^*$ and $v_0'\in\mathbf{V}'$. It follows that
$$
\Psi\Phi\ad_{e_0}=\ad_{\Psi\Phi e_0}\Psi\Phi=\frac1\lambda\ad_{e_0'}\Psi\Phi,
$$
i.e., $(\mathbf{V},\ad_{e_0})\sim(\mathbf{V}',\frac1\lambda\ad_{e_0'})$.
\end{proof}
The finite dimensional analog of this result for $\F\in\{\mathbb{R},\mathbb{C}\}$ is pretty straightforward and can be found in the literature (e.g., \cite{Fre12} or \cite{Gor98}). What we have learnt is that isomorphism classes of almost Abelian Lie algebras correspond to the similarity classes of linear operators on vector spaces $(\mathbf{V},\mbox{T})$. A similarity transformation is merely a change of basis. To describe almost Abelian Lie algebras in isomorphism invariant or basis invariant way we need to find a sufficient similarity invariant for linear operators. Different arts of canonical forms can serve as such invariants. A more explicit description of almost Abelian Lie algebras using canonical forms will be the subject of an upcoming publication.


\section{Derivations and Lie orthogonal operators\label{DerOrth}}

The present section is devoted to derivations and Lie orthogonal operators of almost Abelian Lie algebras.

\begin{proposition}\label{DerBlockFormProp} The algebra of derivations of an almost Abelian Lie algebra $\mathbf{L}=\mathbf{L}_0\oplus\mathbf{W}$ with $\mathbf{L}_0$ indecomposable has the block form
$$
\Der(\mathbf{L})=\left\{\begin{pmatrix}
\phi_{00} & \phi_{01}\\
\phi_{10} & \phi_{11}
\end{pmatrix}\in\End_\F(\mathbf{L})\vline\quad\phi_{00}\in\Der(\mathbf{L}_0),\quad\phi_{01}\mathbf{W}\subset\mathcal{Z}(\mathbf{L}_0),\quad[\mathbf{L}_0,\mathbf{L}_0]\subset\ker\phi_{10}\right\}.
$$
\end{proposition}
\begin{proof} Take $\forall X,X'\in\mathbf{L}_0$ and $Y,Y'\in\mathbf{W}$, and simply write the derivation condition (the Leibnitz rule)
$$
\phi[X+Y,X'+Y']=\phi[X,X']=\phi_{00}[X,X']+\phi_{10}[X,X']=
$$
$$
[\phi(X+Y),X'+Y']+[X+Y,\phi(X'+Y')]=[\phi_{00}X+\phi_{01}Y,X']+[X,\phi_{00}X'+\phi_{01}Y'].
$$
As $X,X',Y,Y'$ are arbitrary we establish the desired result by collecting similar terms.
\end{proof}
\begin{lemma}\label{betaneq0HeisLemma} Let $\mathbf{L}=\mathbf{V}\rtimes\F e_0$ be an indecomposable almost Abelian Lie algebra. If there exists a nonzero linear functional $\beta^\top\in\mathbf{V}^*$ such that $[\mathbf{L},\mathbf{L}]\subset\ker\beta^\top$ and
$$
(\beta^\top v)[e_0,v']=(\beta^\top v')[e_0,v],\quad\forall v,v'\in\mathbf{V}
$$
then $\mathbf{L}_0\simeq\mathbf{H}_\F$.
\end{lemma}
\begin{proof} If $\beta^\top\neq0$ then $\exists v_0\in\mathbf{V}$ such that $\beta^\top v_0\neq0$. From $[\mathbf{L},\mathbf{L}]\subset\ker\beta^\top$ we know that $v_0\notin[\mathbf{L},\mathbf{L}]$. That $\mathbf{L}$ is indecomposable implies by Remark \ref{LDecompRemark} that $\mathcal{Z}(\mathbf{L})\subset[\mathbf{L},\mathbf{L}]$ and therefore $v_0\notin\mathcal{Z}(\mathbf{L})$, i.e., $[e_0,v_0]\neq0$. Denoting $v_1\doteq[e_0,v_0]$ we use the hypothesis of the lemma to write
$$
[e_0,\mathbf{V}]=[\mathbf{L},\mathbf{L}]=\frac{[e_0,v_0]}{\beta^\top v_0}\beta^\top\mathbf{V}=\F v_1
$$
and
$$
[\mathbf{L},\mathbf{L}]\subset\ker\beta^\top=\ker\ad_{e_0}=\mathcal{Z}(\mathbf{L})\subset[\mathbf{L},\mathbf{L}],
$$
whence
$$
[\mathbf{L},\mathbf{L}]=\mathcal{Z}(\mathbf{L})=\F v_1.
$$
But $\ker\beta^\top\subset\mathbf{V}$ is a codimension 1 subspace, so that
$$
\mathbf{V}=\F v_0\oplus_\F\ker\beta^\top=\F v_0\oplus_\F\mathcal{Z}(\mathbf{L})=\F v_0\oplus_\F\F v_1.
$$
Now it is clear that $\mathbf{L}=(\F v_0\oplus\F v_1)\rtimes\F e_0$ is isomorphic to $\mathbf{H}_\F$ (see the final part of the proof of Proposition \ref{2cod1AbIdProp}).
\end{proof}
\begin{proposition} The algebra of derivations of the Heisenberg algebra in the basis of Example \ref{HeisExamp} is
$$
\Der(\mathbf{H}_\F)=\left\{\begin{pmatrix}
\alpha+\Delta_{22} & \Delta_{12} & \gamma_1\\
0 & \Delta_{22} & \gamma_2\\
0 & \beta_2 & \alpha
\end{pmatrix}\vline\quad\alpha,\beta_2,\gamma_1,\gamma_2,\Delta_{12},\Delta_{22}\in\F\right\}.
$$
\end{proposition}
\begin{proof} Simply plug in the block form (\ref{EndBlockForm}) into the derivation condition $\phi[X,Y]=[\phi X,Y]+[X,\phi Y]$ for all $X,Y\in\mathbf{H}_\F$.
\end{proof}
Note the similarity between $\Aut(\mathbf{H}_\F)$ and $\Der(\mathbf{H}_\F)$: in the component $\Delta_{11}$ the determinant of the corresponding algebraic minor is replaced by its trace. In general, when $\Aut(\mathbf{H}_\F)$ is a Lie group then $\Der(\mathbf{H}_\F)$ is related to its Lie algebra.

\begin{proposition}\label{DerLProp} The algebra of derivations of an indecomposable almost Abelian Lie algebra $\mathbf{L}=\mathbf{V}\rtimes\F e_0$ is
$$
\Der(\mathbf{L})=\left\{\begin{pmatrix}
\Delta & \gamma\\
\beta^\top & \alpha
\end{pmatrix}\vline\quad\begin{array}{ll}
\alpha\in\F,\quad\beta^\top\in\mathbf{V}^*,& \gamma\in\mathbf{V},\quad\Delta\in\End_\F(\mathbf{V}),\\
\quad[\mathbf{L},\mathbf{L}]\subset\ker\beta^\top,& (\Delta-\alpha\id)\ad_{e_0}-\ad_{e_0}\Delta=0
\end{array}\right\}.
$$
If $\mathbf{L}\not\simeq\mathbf{H}_\F$ then $\beta^\top=0$.
\end{proposition}
\begin{proof} We simply substitute the block form (\ref{EndBlockForm}) into the derivation condition $\phi[X,Y]=[\phi X,Y]+[X,\phi Y]$ for all $X,Y\in\mathbf{H}_\F$. Collecting similar terms we get the following three equations,
$$
\beta^\top\ad_{e_0}=0,
$$
$$
(\beta^\top v)\ad_{e_0}v'-(\beta^\top v')\ad_{e_0}v=0,
$$
$$
(\Delta-\alpha\id)\ad_{e_0}-\ad_{e_0}\Delta=0.
$$
The first equation simply means that $\ad_{e_0}\mathbf{V}=[\mathbf{L},\mathbf{L}]\subset\ker\beta^\top$. It only remains to use Lemma \ref{betaneq0HeisLemma} to conclude that if $\mathbf{L}\not\simeq\mathbf{H}_\F$ then $\beta^\top=0$.
\end{proof}
\begin{corollary}\label{OuterDerCorr} If all derivations of an almost Abelian Lie algebra are inner, i.e., $\Der(\mathbf{L})=\ad_\mathbf{L}$, then $\mathbf{L}$ is isomorphic to $\mathbf{ax\!+\!b}_\F$.
\end{corollary}
\begin{proof} If $\mathbf{L}=\mathbf{L}_0\oplus\mathbf{W}$ with $\mathbf{L}_0$ indecomposable as usual, then
$$
\ad_\mathbf{L}=\ad_{\mathbf{L}_0}\oplus\ad_{\mathbf{W}}=\ad_{\mathbf{L}_0}\oplus0,
$$
that is, the operator $\ad_{X+Y}$ in (\ref{adRep}) is the block sum of $\ad_X$ and zero for every $X\in\mathbf{L}_0$, $Y\in\mathbf{W}$. But from Proposition \ref{DerBlockFormProp} we know that if $\mathbf{W}\neq0$ then there always exist non-trivial derivations with $\phi_{00}=\phi_{01}=\phi_{10}=0$ and $\phi_{11}\neq0$. Therefore if $\Der(\mathbf{L})=\ad_\mathbf{L}$ then $\mathbf{L}$ is necessarily indecomposable. Then by Proposition \ref{DerLProp} we have
$$
\begin{pmatrix}
\F\ad_{e_0} & \ad_{e_0}\mathbf{V}\\
0 & 0
\end{pmatrix}=\left\{\begin{pmatrix}
\Delta & \gamma\\
\beta^\top & \alpha
\end{pmatrix}\vline\quad\begin{array}{ll}
\alpha\in\F,\quad\beta^\top\in\mathbf{V}^*,& \gamma\in\mathbf{V},\quad\Delta\in\End_\F(\mathbf{V}),\\
\quad[\mathbf{L},\mathbf{L}]\subset\ker\beta^\top,& (\Delta-\alpha\id)\ad_{e_0}-\ad_{e_0}\Delta=0
\end{array}\right\}.
$$
In particular, for $\alpha=0$ and $\beta^\top=0$ the set $\Delta\in\F\ad_{e_0}$ is the full set of solutions of the equation $[\Delta,\ad_{e_0}]=0$. It is clear that $\F\id$ commutes with every operator and hence $\F\id\subset\F\ad_{e_0}$ which means $\ad_{e_0}\in\F^*\id$. But then every $\Delta\in\End_\F(\mathbf{V})$ commutes with $\ad_{e_0}$, hence $\F\id=\End_\F(\mathbf{V})$. This is possible only when $\dim_\F\mathbf{V}=1$, so that $\mathbf{V}=\F v$ for some $v\neq0$ such that $[e_0,v]=\lambda v$ with $\lambda\in\F^*$. It is now not difficult to check directly that the map
$$
(ae_0,bv)\mapsto\begin{pmatrix}
0 & 0\\
b & \lambda a
\end{pmatrix},\quad\forall a,b\in\F
$$
provides a Lie algebra isomorphism between $\mathbf{L}$ and $\mathbf{ax\!+\!b}_\F$.
\end{proof}

Let us now switch to Lie orthogonal operators. The finite dimensional analogue of the result below can be found in \cite{Pop13}, but there the authors use matrix terminology which is not directly applicable in our generality. Recall that $\mathrm{SL}(2,\F)$ is the group of $2\times2$ $\F$-valued matrices of determinant $1$.
\begin{proposition}\label{LieOrthax+bHeisProp} The Lie orthogonal operators of the Lie algebras $\mathbf{ax\!+\!b}_\F$ and $\mathbf{H}_\F$ are
$$
\LO(\mathbf{ax\!+\!b}_\F)=\mathrm{SL}(2,\F),
$$
$$
\LO(\mathbf{H}_\F)=\left\{\begin{pmatrix}
\Delta_{11} & \Delta_{12} & \gamma_1\\
0 & \Delta_{22} & \gamma_2\\
0 & \beta_2 & \alpha
\end{pmatrix}\vline\quad\alpha,\beta_2,\gamma_1,\gamma_2,\Delta_{11},\Delta_{12},\Delta_{22}\in\F,\quad\alpha\Delta_{22}-\gamma_2\beta_2=1\right\}.
$$
\end{proposition}
\begin{proof} These results can be obtained by a direct substitution of the matrix form (\ref{EndBlockForm}) into the Lie orthogonality condition $[\phi X,\phi Y]=[X,Y]$ for all $X$ and $Y$.
\end{proof}

\begin{lemma}\label{1dimadVLemma} For an indecomposable almost Abelian Lie algebra $\mathbf{L}$ if $\dim_\F[\mathbf{L},\mathbf{L}]=1$ then either $\mathbf{L}\simeq\mathbf{ax\!+\!b}_\F$ or $\mathbf{L}\simeq\mathbf{H}_\F$.
\end{lemma}
\begin{proof} Let $\mathbf{L}=\mathbf{V}\rtimes\F e_0$ and let $v_0\in\mathbf{V}$ such that $[\mathbf{L},\mathbf{L}]=\F v_0$. There exists a nonzero $v_1\in\mathbf{V}$ such that $[e_0,v_1]=v_0$. For an arbitrary $v\in\mathbf{V}$, $[e_0,v]=\mu(v)v_0$ for a linear functional $\mu\in\mathbf{V}^*$. It follows that $[e_0,v-\mu(v)v_1]=0$ which means that $v-\mu(v)v_1\in\mathcal{Z}(\mathbf{L})$ for all $v\in\mathbf{V}$. This in turn implies
\begin{equation}
\mathbf{V}=\F v_1\oplus_\F\mathcal{Z}(\mathbf{L}).\label{1dimadVDecomp}
\end{equation}
Because $\mathbf{L}$ is indecomposable we have $\mathcal{Z}(\mathbf{L})\subset\F v_0$. We have two possibilities. If $\mathcal{Z}(\mathbf{L})=0$ then $v_0\notin\mathcal{Z}(\mathbf{L})$ and therefore $[e_0,v_0]=\lambda v_0$ for some $\lambda\in\F^*$, i.e., $v_1=\frac1\lambda v_0$. By (\ref{1dimadVDecomp}) we get $\mathbf{L}=\F v_0\rtimes\F e_0$ which is isomorphic to $\mathbf{ax\!+\!b}_\F$ (compare with the final part of the proof of Corollary \ref{OuterDerCorr}). Otherwise, if $\mathcal{Z}(\mathbf{L})=\F v_0$ then again by (\ref{1dimadVDecomp}) we have $\mathbf{L}=(\F v_0\oplus_\F\F v_1)\rtimes\F e_0$ which is isomorphic to $\mathbf{H}_\F$ (compare with the final part of the proof of Proposition \ref{2cod1AbIdProp}).
\end{proof}

\begin{proposition}\label{LieOrthIndecompLNotHeissNotax+bProp} If an indecomposable almost Abelian Lie algebra $\mathbf{L}$ is not isomorphic to either $\mathbf{ax\!+\!b}_\F$ or $\mathbf{H}_\F$ then its Lie orthogonal operators are
$$
\LO(\mathbf{L})=\left\{\begin{pmatrix}
\frac1\alpha\id+\operatorname{Z} & \gamma\\
0 & \alpha
\end{pmatrix}\vline\quad\alpha\in\F^*,\quad\gamma\in\mathbf{V},\quad\operatorname{Z}\in\Hom(\mathbf{V},\ker\ad_{e_0})\right\}.
$$
\end{proposition}
\begin{proof} As usual, we start by substituting the matrix form (\ref{EndBlockForm}) into the Lie orthogonality condition,
$$
\left[\begin{pmatrix}
\Delta & \gamma\\
\beta^\top & \alpha
\end{pmatrix}\begin{pmatrix}
v\\
t
\end{pmatrix},\begin{pmatrix}
\Delta & \gamma\\
\beta^\top & \alpha
\end{pmatrix}\begin{pmatrix}
v'\\
t'
\end{pmatrix}\right]=\left[\begin{pmatrix}
v\\
t
\end{pmatrix},\begin{pmatrix}
v'\\
t'
\end{pmatrix}\right]=\begin{pmatrix}
t\ad_{e_0}v'-t'\ad_{e_0}v\\
0
\end{pmatrix}.
$$
This is equivalent to the following system of two equations,
\begin{equation}
(\beta^\top v)\ad_{e_0}\Delta v'=(\beta^\top v')\ad_{e_0}\Delta v,\label{LieOrthEq1}
\end{equation}
\begin{equation}
\alpha\ad_{e_0}\Delta v=(\beta^\top v)\ad_{e_0}\gamma+\ad_{e_0}v,\quad\forall v,v'\in\mathbf{V}.\label{LieOrthEq2}
\end{equation}
Suppose first that $\beta^\top\neq0$, which means that there exists $v_0\in\mathbf{V}$ such that $\beta^\top v_0\neq 0$. From (\ref{LieOrthEq1}) we get
$$
\ad_{e_0}\Delta=\frac{\ad_{e_0}\Delta v_0}{\beta^\top v_0}\beta^\top,
$$
and combining with (\ref{LieOrthEq2}) we find that
$$
\ad_{e_0}=\left(\frac{\alpha\ad_{e_0}\Delta v_0}{\beta^\top v_0}-\ad_{e_0}\gamma\right)\beta^\top.
$$
It follows that $[\mathbf{L},\mathbf{L}]=\ad_{e_0}\mathbf{V}$ is 1-dimensional ($\mathbf{L}$ is non-Abelian hence $\ad_{e_0}\neq0$). According to Lemma \ref{1dimadVLemma} we conclude that $\mathbf{L}$ is either $\mathbf{ax\!+\!b}_\F$ or $\mathbf{H}_\F$ which contradicts the synopsis of this proposition. Therefore we necessarily have $\beta^\top=0$. Now (\ref{LieOrthEq1}) becomes trivial and (\ref{LieOrthEq2}) simplifies to
$$
\ad_{e_0}=\alpha\ad_{e_0}\Delta.
$$
This is equivalent to saying that $\alpha\neq 0$ and $\ad_{e_0}(\Delta-\frac1\alpha\id)=0$. Denoting $\mathrm{Z}\doteq\Delta-\frac1\alpha\id$ we see that $\mathrm{Z}\mathbf{V}\subset\ker\ad_{e_0}=\mathcal{Z}(\mathbf{L})$ which completes the proof.
\end{proof}

\begin{corollary}\label{phiL+ZL=LCorr} If $\mathbf{L}$ is an indecomposable almost Abelian Lie algebra and if $\phi\in\LO(\mathbf{L})$ is any Lie orthogonal operator then $\phi\mathbf{L}+\mathcal{Z}(\mathbf{L})=\mathbf{L}$.
\end{corollary}
\begin{proof} By Proposition \ref{LieOrthax+bHeisProp} the claim is obvious for $\mathbf{ax\!+\!b}_\F$ as all Lie orthogonal operators are bijective. For $\mathbf{H}_\F$ a Lie orthogonal operator $\phi$ in the form presented in Proposition \ref{LieOrthax+bHeisProp} is non-surjective when $\Delta_{11}=0$, but it is precisely the centre $\mathcal{Z}(\mathbf{H}_\F)$ that misses in the image $\phi\mathbf{H}_\F$ and therefore again $\mathbf{H}_\F=\phi\mathbf{H}_\F+\mathcal{Z}(\mathbf{H}_\F)$. Now assume $\mathbf{L}$ is not one of the two algebras above. Then by Proposition \ref{LieOrthIndecompLNotHeissNotax+bProp} we have that
$$
\phi\mathbf{L}=\left(\F\gamma+\left(\frac1\alpha\id+\mbox{Z}\right)\mathbf{V}\right)\oplus_\F\F e_0.
$$
The identity
$$
v=\left(\frac1\alpha\id+\mathrm{Z}\right)\alpha v-\alpha\mathrm{Z}v,\quad\forall v\in\mathbf{V}
$$
along with the fact that $\mathrm{Z}\mathbf{V}\in\mathcal{Z}(\mathbf{L})$ show that
$$
\mathbf{V}=\left(\frac1\alpha\id+\mathrm{Z}\right)\mathbf{V}+\mathcal{Z}(\mathbf{L})
$$
which completes the proof.
\end{proof}

\begin{proposition} For an almost Abelian Lie algebra $\mathbf{L}=\mathbf{L}_0\oplus\mathbf{W}$ with $\mathbf{L}_0$ indecomposable the Lie orthogonal operators have the following block form,
$$
\LO(\mathbf{L})=\left\{\begin{pmatrix}
\phi_{00} & \phi_{01}\\
\phi_{10} & \phi_{11}
\end{pmatrix}\in\End_\F(\mathbf{L})\vline\quad\phi_{00}\in\LO(\mathbf{L}_0),\quad\phi_{01}\mathbf{W}\subset\mathcal{Z}(\mathbf{L}_0)\right\}.
$$
\end{proposition}
\begin{proof} We write down the Lie orthogonality condition $[\phi(X+Y),\phi(X'+Y')]=[X+Y,X'+Y']$ for arbitrary $X,X'\in\mathbf{L}_0$ and $Y,Y'\in\mathbf{W}$. Collecting similar terms we get the following pair of equations,
$$
[\phi_{00}X,\phi_{00}X']=[X,X'],\quad[\phi_{01}\mathbf{W},\phi_{00}\mathbf{L}_0+\phi_{01}\mathbf{W}]=0.
$$
The first equation simply means $\phi_{00}\in\LO(\mathbf{L}_0)$. The second equation in view of Corollary \ref{phiL+ZL=LCorr} gives
$$
[\phi_{01}\mathbf{W},\mathbf{L}_0]=[\phi_{01}\mathbf{W},\phi_{00}\mathbf{L}_0+\mathcal{Z}(\mathbf{L}_0)]=[\phi_{01}\mathbf{W},\phi_{00}\mathbf{L}_0]=0
$$
whence we get $\phi_{01}\mathbf{W}\subset\mathcal{Z}(\mathbf{L}_0)$.
\end{proof}


\section{The centre of the universal enveloping algebra\label{Casimir}}

Here we are going to find the centre of the universal enveloping algebra of an almost Abelian Lie algebra $\mathbf{L}$. If $\mathfrak{h}:\mathbf{L}\mapsto\operatorname{U}(\mathbf{L})$ is the embedding map and $\mathbf{L}=\mathbf{V}\rtimes\F e_0$ is an almost Abelian Lie algebra, then let us choose a basis $\{e_i\}_{i\in\aleph}$ of $\mathbf{V}$ and denote
$$
x_0=\mathfrak{h}(e_0),\quad x_i=\mathfrak{h}(e_i),\quad\forall i\in\aleph.
$$
Then $x_ix_j=x_jx_i$ for all $i,j\in\aleph$. By Poincar\'e-Birkhoff-Witt theorem the monomials of the form $x_0^nP(x_\aleph)$ with $P(x_\aleph)\in\F[\{x_i\}_{i\in\aleph}]$ span $\operatorname{U}(\mathbf{L})$. In other words, every element $x\in\operatorname{U}(\mathbf{L})$ has a unique representation
$$
x=\sum_{m=0}^nx_0^mP_m(x_\aleph),\quad n\in\mathbf{N}_0,\quad P_m(x_\aleph)\in\F[\{x_i\}_{i\in\aleph}].
$$
Denote by $\ad\in\Der(\operatorname{U}(\mathbf{L}))$ the extension of the derivation $\ad_{e_0}\in\Der(\mathbf{L})$. For $\ad$ the variable $x_0$ is a constant,
$$
\ad x=\sum_{m=0}^nx_0^m\ad P_m(x_\aleph),
$$
and on every $P(x_\aleph)$ is acts as a derivation,
$$
\ad (x_{i_1}x_{i_2}...x_{i_p})=(\ad x_{i_1})x_{i_2}...x_{i_p}+x_{i_1}(\ad x_{i_2})...x_{i_p}+...+x_{i_1}x_{i_2}...(\ad x_{i_p}).
$$
In particular, $\ad$ preserves the degree of every homogeneous polynomial (zero monomial is of any degree of homogeneity).

\begin{lemma}\label{Pxx0CommLemma} For every $n\in\mathbf{N}_0$ and $P(x_\aleph)\in\F[\{x_i\}_{i\in\aleph}]$ we have
$$
P(x_\aleph)x_0^n=(x_0-\ad)^nP(x_\aleph).
$$
\end{lemma}
\begin{proof} We start by noting that
$$
x_ix_0=x_0x_i-\ad x_i,\quad\forall i\in\aleph.
$$
Assume that for every monomial of degree $q$ less than some $p\in\mathbb{N}$ it holds
$$
x_{i_1}...x_{i_q}x_0=(x_0-\ad)(x_{i_1}...x_{i_q}).
$$
Then for every monomial of degree $p$ we have
$$
x_{i_1}...x_{i_p}x_0=x_{i_1}(x_0-\ad)(x_{i_2}...x_{i_p})=x_{i_1}x_0x_{i_2}...x_{i_p}-x_{i_1}\ad(x_{i_2}...x_{i_p})=
$$
$$
(x_0-\ad)(x_{i_1})x_{i_2}...x_{i_p}-x_{i_1}\ad(x_{i_2}...x_{i_p})=(x_0-\ad)x_{i_1}...x_{i_p}.
$$
Thus, by mathematical induction we prove the statement for all $p\in\mathbb{N}$. Every polynomial is a sum of monomials, hence
$$
P(x_\aleph)x_0=(x_0-\ad)P(x_\aleph),\quad\forall P(x_\aleph)\in\F[\{x_i\}_{i\in\aleph}].
$$
Now assume that
$$
P(x_\aleph)x_0^m=(x_0-\ad)^mP(x_\aleph)
$$
holds for all $P(x_\aleph)$ and $m<n$ for a fixed $n\in\mathbb{N}$. Then
$$
P(x_\aleph)x_0^n=P(x_\aleph)x_0^{n-1}x_0=(x_0-\ad)^{n-1}P(x_\aleph)x_0=(x_0-\ad)^nP(x_\aleph).
$$
Using mathematical induction once more we arrive at the desired result.
\end{proof}

\begin{proposition} The centre of the universal enveloping algebra of an almost Abelian Lie algebra $\mathbf{L}$ as above consists of $\ad$-conserved polynomials in $\{x_i\}_{i\in\aleph}$,
$$
\mathcal{Z}(\operatorname{U}(\mathbf{L}))=\left\{P(x_\aleph)\in\F[\{x_i\}_{i\in\aleph}]|\quad\ad P(x_\aleph)=0\right\}.
$$
\end{proposition}
\begin{proof} Let $x\in\operatorname{U}(\mathbf{L})$ be an arbitrary element of the universal enveloping algebra,
$$
x=x_0^nP_n(x_\aleph)+...+P_0(x_\aleph),\quad n\in\mathbf{N}_0,\quad P_m(x_\aleph)\in\F[\{x_i\}_{i\in\aleph}].
$$
The condition that $x\in\mathcal{Z}(\operatorname{U}(\mathbf{L}))$ is equivalent to the following set of two conditions,
$$
xx_0=x_0x,\quad xx_i=x_ix,\quad\forall i\in\aleph.
$$
Using Lemma \ref{Pxx0CommLemma} the first condition can be transformed as
$$
xx_0=(x_0^nP_n(x_\aleph)+...+P_0(x_\aleph))x_0=x_0^n(x_0-\ad)P_n(x_\aleph)+...+(x_0-\ad)P_0(x_\aleph)=
$$
$$
x_0^{n+1}P_n(x_\aleph)+...+x_0P_0(x_\aleph)-x_0^n\ad P_n(x_\aleph)-...-\ad P_0(x_\aleph)
$$
$$
=x_0x-x_0^n\ad P_n(x_\aleph)-...-\ad P_0(x_\aleph)=x_0x,
$$
which is equivalent to
$$
\ad P_m(x_\aleph)=0,\quad m=0,...,n.
$$
The second condition becomes (again with the help of Lemma \ref{Pxx0CommLemma})
$$
x_ix=x_i(x_0^nP_n(x_\aleph)+...+P_0(x_\aleph))=[(x_0-\ad)^nx_i]P_n(x_\aleph)+...+x_iP_0(x_\aleph)=
$$
$$
\left[x_0^nx_i+\sum_{l=1}^nx_0^{n-l}(-\ad)^lx_i\right]P_n(x_\aleph)+...+x_iP_0(x_\aleph)=
$$
$$
x_0^nP_n(x_\aleph)x_i+...+P_0(x_\aleph)x_i+\sum_{l=1}^nx_0^{n-l}\left[(-\ad)^lx_i\right]P_n(x_\aleph)+...+\sum_{l=1}^1x_0^{1-l}\left[(-\ad)^lx_i\right]P_1(x_\aleph)=
$$
$$
xx_i+\sum_{l=1}^nx_0^{n-l}\left[(-\ad)^lx_i\right]P_n(x_\aleph)+...+\sum_{l=1}^1x_0^{1-l}\left[(-\ad)^lx_i\right]P_1(x_\aleph)=xx_i,\quad\forall i\in\aleph,
$$
which means
$$
\sum_{l=1}^nx_0^{n-l}\left[(-\ad)^lx_i\right]P_n(x_\aleph)+...+\sum_{l=1}^1x_0^{1-l}\left[(-\ad)^lx_i\right]P_1(x_\aleph)=
$$
$$
x_0^{n-1}\left[(-\ad)x_i\right]P_n(x_\aleph)+x_0^{n-2}\left(\left[(-\ad)^2x_i\right]P_n(x_\aleph)+\left[(-\ad)x_i\right]P_{n-1}(x_\aleph)\right)+...
$$
$$
+\left[(-\ad)^nx_i\right]P_n(x_i)+...+\left[(-\ad)x_i\right]P_1(x_\aleph)=0,\quad\forall i\in\aleph.
$$
Start with the term proportional to $x_0^{n-1}$. Because $\mathbf{L}$ is non-Abelian, there exists $i\in\aleph$ such that $\ad x_i\neq0$ and thus $P_n(x_\aleph)=0$. Take next the term proportional to $x_0^{n-2}$ and set $P_n(x_\aleph)=0$: we get $P_{n-1}(x_\aleph)=0$. Continuing in a similar fashion we will find that
$$
P_m(x_\aleph)=0,\quad m=1,...,n.
$$
Therefore
$$
x=P_0(x_\aleph),\quad\ad P_0(x_\aleph)=0
$$
as asserted.
\end{proof}

\begin{remark} Note that if $X\in\mathcal{Z}(\mathbf{L})=\ker\ad_{e_0}$ and therefore $x=\mathfrak{h}(X)\in\ker\ad$ then $\F[x]\in\mathcal{Z}(\operatorname{U}(\mathbf{L}))$. We can write this symbolically as $\F[\mathcal{Z}(\mathbf{L})]\subset\mathcal{Z}(\operatorname{U}(\mathbf{L}))$. The next example shows that this inclusion is often proper.
\end{remark}

\begin{example} Let $\mbox{\normalfont Bi}(\mbox{\normalfont VII}_0)=\F^2\rtimes\F e_0$ be the Bianchi $\mbox{\normalfont VII}_0$ Lie algebra,
$$
\mbox{\normalfont Bi}(\mbox{\normalfont VII}_0)=\left\{\begin{pmatrix}
0 & 0 & 0\\
v & 0 & -t\\
u & t & 0
\end{pmatrix}\vline\quad(t,v,u)\in\F^3\right\}.
$$
This corresponds to
$$
\ad_{e_0}=\begin{pmatrix}
0 & -1\\
1 & 0
\end{pmatrix}.
$$
If we let $\operatorname{U}(\mbox{\normalfont Bi}(\mbox{\normalfont VII}_0))=\F[x_0,x_1,x_2]$ with
$$
x_0=\mathfrak{h}\begin{pmatrix}
1\\
0\\
0
\end{pmatrix},\quad x_1=\mathfrak{h}\begin{pmatrix}
0\\
1\\
0
\end{pmatrix},\quad x_2=\mathfrak{h}\begin{pmatrix}
0\\
0\\
1
\end{pmatrix},
$$
then we have $\ad x_1=x_2$ and $\ad x_2=-x_1$. It follows that $\ad(x_1^2+x_2^2)=0$, i.e., $x_1^2+x_2^2\in\mathcal{Z}(\operatorname{U}(\mbox{\normalfont Bi}(\mbox{\normalfont VII}_0)))$ although $\mathcal{Z}(\mbox{\normalfont Bi}(\mbox{\normalfont VII}_0))=\ker\ad_{e_0}=0$.
\end{example}


\end{document}